\documentclass{article}
\usepackage[cp1251]{inputenc}
\usepackage[english]{babel}
\usepackage[12pt]{extsizes}
\usepackage{mathrsfs}
\usepackage{amsfonts,amssymb}
\usepackage[tbtags]{amsmath}
\usepackage[T2A]{fontenc}
\usepackage[unicode, pdftex]{hyperref}
\usepackage{mathtools}
\usepackage{blindtext,titlefoot}

\DeclarePairedDelimiter{\floor}{\lfloor}{\rfloor}

\DeclarePairedDelimiter{\ceil}{\lceil}{\rceil}

\newtheorem{theorem}{Theorem}
\newtheorem{definition}{Definition}
\newtheorem{lemma}{Lemma}
\newtheorem{claim}{Claim}

\begin{document}
\title{Tight asymptotics of clique-chromatic numbers of dense random graphs}
\author{Yu. Demidovich\footnote{King Abdullah University of Science and Technology; KAUST Artificial Intelligence Initiative, e-mail: yury.demidovich@kaust.edu.sa} \and M. Zhukovskii\footnote{University of Sheffield, Department of Computer Science, e-mail: zhukmax@gmail.com}}
\date{}
\maketitle
\begin{abstract}
	The clique chromatic number of a graph is the minimum number of colors required to assign to its vertex set so that no inclusion maximal clique is monochromatic. McDiarmid, Mitsche and Pra\l at proved that the clique chromatic number of the binomial random graph $G\left(n,\frac{1}{2}\right) $ is at most $\left(\frac{1}{2}+o(1)\right)\log_2n$ with high probability. Alon and Krivelevich showed that it is greater than $\frac{1}{2000}\log_2n$ with high probability and suggested that the right constant in front of the logarithm is $\frac{1}{2}.$ We prove their conjecture and, beyond that, obtain a tight concentration result: whp
	$\chi_c\left(G\left(n,\frac{1}{2}\right)\right) = \frac{1}{2}\log_2 n - \Theta\left(\ln\ln n\right).$
\end{abstract}
\section{Introduction}
\textit{A clique coloring} of a graph $G=(V,E)$ is a coloring of $V$ such that there is no inclusion maximal monochromatic clique (i.e. every monochromatic clique has a vertex adjacent to all its vertices and colored in another color). \textit{The clique chromatic number} $\chi_c(G)$ is the minimum number of colors in a clique coloring of $G.$

Clearly, $\chi_c(G)$ is at most the ordinary chromatic number of $G$ (i.e. the minimum number of colors in a coloring without monochromatic edges) denoted by $\chi(G).$ Note, that if $G$ is triangle-free, then $\chi_c(G) = \chi(G).$ However, the clique-chromatic number of a graph can be much smaller than its chromatic number. For instance, $\chi_c(K_{n}) = 2$ while $\chi(K_n) = n.$ Futhermore, some decision problems formulated in terms of the chromatic number are harder for their reformulations in terms of the clique chromatic number. In particular, deciding whether $\chi_c(G)=2$ is NP-complete (see \cite{BGG, KrT}) while determining whether or not $G$ is bipartite (i.e. $\chi(G)=2$) is computable in linear time. Another important difference between these characteristics is that unlike the usual chromatic number, the clique chromatic number is not monotonic with respect to taking subgraphs. Many additional structural results on the clique chromatic number can be found in \cite{AST, CPTT,JMRS,MohS}. A fundamental question about a graph characteristic is its behavior on (asymptotically) almost all graphs. It is known (see \cite{Bol}) that almost all graphs have chromatic number $\frac{n}{2\log n}(1+o(1))$ (everywhere below we write $\log n$ for the logarithm to the base of $2$ and $\ln n$ for the natural logarithm). Recently in \cite{HecRio} the long standing question about the concentration interval for the chromatic number was almost solved: roughly, the optimal size of the concentration interval for the chromatic number of almost all graphs is $n^{1/2-o(1)}.$ In this paper, we find the clique chromatic number of almost all graphs and give a tight bound for the size of the concentration interval.$\\$

Let us recall that \textit{the random graph} $G\left(n,\frac{1}{2}\right)$ is a random uniformly distributed element of the set of all graphs on $[n]:=\{1,\ldots,n\}$ (or, in other words, every pair of vertices is adjacent with probability $p=\frac{1}{2}$ independently). We say that a graph property $Q$ holds with high probability (whp), if ${\sf P}\left(G\left(n,\frac{1}{2}\right)\in Q \right) \to 1$ as $n\to\infty.$

In 2016 \cite{MMP1}, McDiarmid, Mitsche and Pra\l at proved that whp
\begin{equation}\label{upper_MMP}
	\chi_c\left(G\left(n,\frac{1}{2}\right)\right)\leqslant\ceil[\bigg]{\left(\frac{1}{2}+\frac{2\log\ln n}{\log n}\right)\log n}+1
\end{equation}  
and asked about a lower bound. In 2017 \cite{AK}, Alon and Krivelevich proved that $\log n$ is the right order of magnitude by showing that whp  $\chi_c\left( G\left(n,\frac{1}{2}\right)\right) \geqslant\frac{1}{2000}\log n.$ They also suggested that the right constant in front of $\log n$ is $\frac{1}{2}.$ In the present work we refine the upper bound \eqref{upper_MMP} and show that the conjecture of Alon and Krivelevich holds true.
\begin{theorem}\label{MainResult}
	With high probability
	\begin{equation*}
		\frac{1}{2}\log n - \left(3+o(1)\right)\log\ln n \leqslant\chi_c\left(G\left(n,\frac{1}{2}\right)\right)\leqslant \frac{1}{2}\log n - \left(\frac{1}{2}-o(1)\right) \log\ln n.
	\end{equation*} 
\end{theorem}

The upper bound is easy. It is simply obtained by the naive (greedy) algorithm proposed in \cite[Theorem 3.4]{MMP1}: at step $i$, assign the new color $i$ to all uncolored neighbors of $i,$ and proceed in this way until the set of all uncolored vertices does not contain a maximal clique. We estimate the number of steps more carefully, and thus obtain a slightly better bound than \eqref{upper_MMP}. The proof is given in Section \ref{section_upper_bound}.

Section \ref{section_lower_bound} is devoted to the proof of the lower bound which is the trickiest part of the paper. On the one hand, the proof is (unsurprisingly) based on the following natural property of almost all graphs which is claimed in Lemma \ref{lemma}, Section \ref{lemma_key}: if we fix a large enough set $Y\subset[n],$ then whp, if all the vertices outside this set have many non-neighbors inside it, then $Y$ contains a maximal clique. On the other hand, the probability of this property tends to $1$ not fast enough, so the union bound (over $Y$) does not work. However, given a coloring in $s$ colors of a graph $G$ for small enough $s$, we propose an algorithm for selecting a subset of vertices of $G$ that can not be colored in the desired number of colors (see Section \ref{final_step}), if the graph $G$ satisfies a certain property presented in Definition \ref{def1} in Section \ref{property}. This property describes a family of large sets $Y\subset[n]$ with the above mentioned requirement: it might contain a maximal clique if every vertex outside $Y$ has many non-neighbors inside. The cardinality of this family is small enough, and this allows to apply the union bound over $Y$ in this family together with the probability bound from Lemma \ref{lemma}. It is also worth mentioning that the presented algorithm of selecting an uncolored set implies that all clique colorings in not too many colors are close to the ``greedy coloring'' (see Section \ref{conclusion}). Notice that a weaker version of the mentioned property of a subset $Y$ was used by Alon and Krivelevich in \cite{AK}. However, their version is not sufficient (even after choosing an optimal set of parameters and after applying our algorithm), see the outline in Section \ref{key_lemma}.

The clique chromatic number of sparse binomial random graphs (i.e. $p\to 0$ as $n\to\infty$) was also studied. Asymptotical behavior was tightly estimated in \cite{MMP1} for most values of $p,$ and these results were refined in \cite{LMW}.
$\\$
 
Before switching to the proofs, let us introduce several notations that we use in the next sections. Let $N_0(v_1,\ldots,v_j)$ be the set of common non-neighbors of $v_1,\ldots,v_j$ (the host graph is always clear from the context). For a graph $G,$ denote by $G\vert_{U}$ the subgraph of $G$ induced on $U\subseteq V(G).$

\section{Upper bound}\label{section_upper_bound}

Let $\varepsilon>0$ be sufficiently small. Let $s = \ceil[\big]{\frac{1}{2}\log n  - \left(\frac{1}{2} - \varepsilon\right)\log\ln n}.$ Take vertices $1,\ldots,s.$ By Chernoff's bound \cite[Theorem 2.1]{JLR}, whp
\begin{equation*}
	\left|N_0\left(1,\ldots,s\right)\right|\leqslant \frac{n}{2^s} + \sqrt{\frac{n}{2^s}}\ln n\leqslant \sqrt{n}\left(\ln n\right)^{\frac{1}{2} - \varepsilon} + \sqrt[4]{n}\left(\ln n\right)^{5/4-\varepsilon/2}.
\end{equation*}
It is clear, that $G(n,\frac{1}{2})\vert_{N_0(1,\ldots,s)}\stackrel{d}=G\left(|N_0(1,\ldots,s)|,\frac{1}{2}\right).$

Let us describe the greedy algorithm. At the first step, assign the color $1$ to all neighbors of $1.$ At every step $i=2,…,s$ assign color $i$ to all uncolored neighbors of $i$. After $s$ steps, the uncolored vertices among $1,…,s$ form an independent set. Assign to them color $s+1.$ It remains to show that there is no maximal clique inside the set of uncolored vertices. For that, we need
\begin{claim}\label{lemma_upper}
	Whp, for every positive integer $m\leqslant\sqrt{n}(\ln n)^{1/2 - \varepsilon}+\sqrt[4]{n}\left(\ln n\right)^{5/4 - \varepsilon / 2},$ $G\left(n,\frac{1}{2}\right)\vert_{\left[m\right]}$ does not contain a maximal clique.
\end{claim}
The proof is given in the Appendix since it is based on the direct application of Markov's inequality and does not give any additional insight.

By Claim \ref{lemma_upper}, whp $N_0(1,\ldots,s)$ does not contain a maximal clique. Therefore, it suffices to assign to every vertex of $N_0(1,\ldots,s)$ only one additional color. It follows that whp
\begin{equation*}
	\chi_c\left(G\left(n,\frac{1}{2}\right)\right)\leqslant s + 2 = \\ \ceil[\bigg]{\frac{1}{2}\log n  - \left(\frac{1}{2} - \varepsilon\right)\log\ln n} + 2.
\end{equation*}

\section{Lower bound}\label{section_lower_bound}
\subsection{The property of Y}\label{lemma_key}
Set $N_0\left(\varnothing\right) = [n].$ Designate by $N_0(v,U)$ the set of non-neighbors of $v$ in $U$ (the host graph is always clear from the context). Fix $\varepsilon>0$ small enough. 
Let 
$$
k = \ceil[\bigg]{\log n + \left(\frac{1}{\ln 2} + 4\varepsilon\right)\ln\ln n}.
$$
\begin{lemma}\label{lemma}
	Fix a set $Y\subseteq [n]$ of size $|Y|:=y\geqslant \left(\ln n\right)^{2+3\varepsilon}\sqrt{n}.$ For $n$ large enough, the probability that
	\begin{itemize}
		\item $|N_0(v,Y)|\geqslant \left(\ln n\right)^{2+2\varepsilon}\sqrt{n}$ for every $v\in [n]\setminus Y,$
		\item there are no cliques $K_{k}\subset Y$ such that every vertex in $[n]\setminus Y$ has a non-neighbor in $K_{k}$
	\end{itemize} is at most $e^{-\left(1+\frac{\varepsilon}{2}\right)y\ln\ln n}.$
\end{lemma}

Lemma \ref{lemma} is proven in Section \ref{key_lemma}. In Section \ref{property} we define a graph property that describes the sufficient set of $Y$'s satisfying the conditions from Lemma \ref{lemma}. This property implies the lower bound in Theorem \ref{MainResult} and allows exploiting the union bound in a smart way so that the probability bound from Lemma \ref{lemma} is sufficient.
In Section \ref{final_step}, we show how this property works for finding a set of uncolored vertices and finish the proof of the lower bound.

\subsection{Threatened sets}\label{property}
Let $\varepsilon>0$ be sufficiently small. Put $s:=\floor[\big]{\frac{1}{2}\log n - \left(\frac{3}{\ln 2} + 5\varepsilon\right)\ln\ln n}.$ Let $G$ be a graph on the vertex set $[n].$
\begin{definition}\label{def1} $G$ has Property \textbf{C}, if, for every $j\in[s]$ and vertices $v_1,\ldots,v_j,$
	\begin{enumerate}
		\item $|N_0(v_1,\ldots,v_j)| \geqslant \frac{n}{2^j}-2\sqrt{n}\ln n,$
		\item any $Y\subseteq N_0(v_1,\ldots,v_j),$ s.t.
		\begin{itemize}
			\item $|Y|\geqslant \frac{|N_0(v_1,\ldots,v_j)|}{\log n},$
			\item $|N_0(v,Y)|\geqslant \left(\ln n\right)^{2+ 2\varepsilon}\sqrt{n}$ for every $v\in [n]\setminus
			Y,$
		\end{itemize}
		contains a maximal clique.
	\end{enumerate}
\end{definition}
\begin{lemma}\label{lemma2}
	Whp $G(n,\frac{1}{2})$ has Property \textbf{C}.
\end{lemma}

\noindent\textbf{Proof of Lemma \ref{lemma2}.} Let $j\in[s],$ $v_1,\ldots,v_j\in [n].$ Denote by $\xi$ the number of vertices in $G\left(n,\frac{1}{2}\right)$ non-adjacent to any of $v_1,\ldots,v_j.$  Denote by $B$ the event that $\xi \geqslant \frac{n}{2^j}-2\sqrt{n}\ln n.$ Since $\xi$ is distributed as $\mathrm{Bin}\left( n-j,\frac{1}{2^j}\right),$ by Chernoff's bound \cite[Theorem 2.1]{JLR},
\begin{equation*}
	{\sf P}\left(\overline{B} \right) \leqslant e^{-2\ln^2n}.
\end{equation*}

Denote by $A$ the event saying that Property $2$ from Definition \ref{def1} holds for $v_1,\ldots,v_j.$ Expose first all the edges having vertices among $v_1,\ldots,v_j.$ Define $G_j:=G\vert_{ [n]\setminus\{v_1,\ldots,v_j\}}$ $\stackrel{d}= G\left(n-j,\frac{1}{2}\right).$  Choose $Y\subset N_0(v_1,\ldots,v_j)\subset V(G_j)$ of size at least $\frac{\xi}{\log n}.$ Since the distribution of edges of $G_j$ does not depend on neighbors of $v_1,\ldots,v_j,$ we may apply Lemma \ref{lemma} to $G_j$ (if a clique in $Y$ is maximal in $G_j,$ then it is also maximal in $G$) and the union bound (over $Y$):
\begin{multline*}
	{\sf P}\left(\overline{A\cap B} \right)= {\sf P}(\overline{A}\cap B)+{\sf P}(\overline{B})\leqslant {\sf P}(\overline{A}|B)+{\sf P}(\overline{B})\leqslant\\\leqslant \max_{N\geqslant \frac{n}{2^j}-2\sqrt{n}\ln n}\sum_{y\geqslant \frac{N}{\log n}}\binom{N}{y}e^{-\left(1+\frac{\varepsilon}{2}\right)y(1+o(1))\ln\ln n}+e^{-2\ln^2n}\leqslant\\
	\leqslant \max_{N\geqslant \frac{n}{2^j}-2\sqrt{n}\ln n}\sum_{y\geqslant \frac{N}{\log n}}e^{y\left(1 + \ln\log n-\left(1+\frac{\varepsilon}{2}\right)\ln\ln n\right)(1+o(1)) }+e^{-2\ln^2n}\leqslant\\
	\leqslant \max_{N\geqslant \frac{n}{2^j}-2\sqrt{n}\ln n}n e^{-\frac{\varepsilon N\ln\ln n}{2\log n}(1+o(1))}+e^{-2\ln^2n}
	=e^{-2(1+o(1))\ln^2n}.
\end{multline*}

By the union bound, the probability that $G\left(n,\frac{1}{2}\right)$ does not have Property \textbf{C} is at most
\begin{equation*}
	\sum_{j=0}^{s} \binom{n}{j} e^{-2(1+o(1))\ln^2n}\leqslant se^{s\ln n-2(1+o(1))\ln^2n} = se^{\left(\frac{1}{2\ln 2} - 2\right)(1+o(1))\ln^2n}\to 0
\end{equation*}
as $n\to\infty.$
\begin{flushright}
	$\blacksquare$
\end{flushright}
\subsection{Uncolored set}\label{final_step}
Let $G$ be a graph on the vertex set $[n]$ with Property $\textbf{C}.$ Suppose there is a proper assignment of $s=\floor[\big]{\frac{\log n}{2} - \left(\frac{3}{\ln 2}+ 5\varepsilon\right)\ln\ln n}$ colors to $[n].$ Due to Lemma \ref{lemma2}, to prove Theorem \ref{MainResult} it is sufficient to get a contradiction with Property \textbf{C}.

There exists a color class $Y_1$ with at least $\frac{n}{\log n}$ vertices. Find a vertex $v_1\in [n]\setminus Y_1$ that has the smallest number of non-neighbors in $Y_1.$ Since the coloring is proper,  $|N_0(v_1,Y_1)|<\left(\ln n\right)^{2+2\varepsilon}\sqrt{n}$ by Property \textbf{C}.

For $j=1,\ldots,s-1$ repeat the following procedure. 

Put $X=N_0(v_1,\ldots,v_j),$ $x=|X|.$ There exists a color class other than $Y_1,\ldots, Y_j$ that has at least $\frac{x}{\log n}$ vertices in $X,$ choose a class that has the biggest number of vertices in $X$ and denote its intersection with $X$ by $Y_{j+1}.$ Find a vertex $v_{j+1}\in [n]\setminus Y_{j+1}$ (it can not be equal to any of $v_1,\ldots,v_j$ since otherwise $N_0(v_{j+1},Y_{j+1}) = Y_{j+1}$) that has the smallest number of non-neighbors in $Y_{j+1}.$ Since the coloring of $G$ is proper, $|N_0(v_{j+1},Y_{j+1})|<\left(\ln n\right)^{2+2\varepsilon}\sqrt{n}$ by Property \textbf{C}.

At the end of the process, by Property \textbf{C}, we obtain the set $N_0(v_1,\ldots,v_s)$ of size $\Omega\left(\sqrt{n}\left(\ln n\right)^{3+5\varepsilon\ln 2}\right)$ colored in $s$ colors with every color class smaller than $\left(\ln n\right)^{2+2\varepsilon}\sqrt{n}.$ It remains to note that $s\left(\ln n\right)^{2+2\varepsilon}\sqrt{n} = o\left(\left(\ln n\right)^{3+5\varepsilon\ln 2}\sqrt{n}\right)$ --- contradiction.

\subsection{Proof of Lemma \ref{lemma}}\label{key_lemma}
\noindent\textbf{Outline of the proof.} In the proof of Lemma 2.3 from \cite{AK}, Alon and Krivelevich select a small set of ``bad'' vertices in $[n]\setminus Y$ that have a small number of non-neighbors in $Y$ (all the other vertices have roughly $|Y|/2$ non-neighbors in $Y$). After that, for every ``bad'' vertex, they choose its own set of its non-neighbors in $Y.$ The rest of $Y$ (denoted by $Y'$) is divided into parts so that every ``good'' vertex has non-neighbors in, roughly, half of each part. Finally, they include in a maximal clique in $Y$ one vertex from every set of the partition. This approach does not work in our case, because the number of ``bad'' vertices is larger than the size of a maximal clique. We divide the set of ``bad'' vertices into two parts. For the first part (these vertices have extremely small number of non-neighbors in $Y$), we apply the same steps as Alon and Krivelevich. Then, we find a partition of $Y'$ such that every ``good'' vertex has non-neighbors in a half of each set, as before, and each vertex from the second part of the set of ``bad'' vertices has non-neighbors in a fraction of each set.$\\$

\noindent\textbf{Proof.} Let 
$$
t_1 = \floor{\ln^2n},\quad
t_2 =\floor[\bigg]{ \left(\frac{1}{\ln 2}+3\varepsilon\right)\ln\ln n}.
$$
Notice that $x\ln x+(1-x)\ln(1-x)\to 0$ as $x\to 0+.$ Choose $\alpha>0$ so small that $|\alpha\ln\alpha + (1-\alpha)\ln(1-\alpha)|<\varepsilon.$ Let 
\begin{align*}
	B_1=&\left\lbrace v\in[n]\setminus Y: \alpha y < \left|N_0(v,Y)\right|<\left(\frac{1}{2}-\frac{\varepsilon}{5}\frac{\ln\ln n}{\ln n}\right)y \right\rbrace,\\
	B_2=&\left\lbrace v\in[n]\setminus Y: |N_0(v, Y)|\leqslant \alpha y \right\rbrace.	
\end{align*}

For $v\in[n]\setminus Y,$ define 
\begin{equation*}
	q_1:={\sf P}\left(\alpha y < \left|N_0(v,Y)\right|<\left(\frac{1}{2}-\frac{\varepsilon}{5}\frac{\ln\ln n}{\ln n}\right)y\right),\qquad
	q_2:={\sf P}\left(\left|N_0(v,Y)\right|\leqslant\alpha y\right).
\end{equation*}
Clearly, $q_1$ and $q_2$ do not depend on the choice of $v.$ By Chernoff's bound \cite[Theorem 2.1]{JLR}, $q_1\leqslant e^{-\frac{\varepsilon^2\ln^2\ln n}{25\ln^2 n}y}.$ Let us also bound $q_2$ from above. For large $n,$ we have
\begin{multline*}
	{\sf P}\left(\mathrm{Bin}\left(y,\frac{1}{2}\right)<\alpha y\right)\leqslant\alpha y\binom{y}{\ceil{\alpha y}}\left(\frac{1}{2}\right)^y =\\ O\left(\sqrt{n}\right)\cdot\left(\frac{1}{\alpha^{\alpha}(1-\alpha)^{1-\alpha}}\right)^y2^{-y} < e^{-(\ln 2 - \varepsilon)y}.
\end{multline*}
Therefore,	$q_2\leqslant e^{-(\ln 2 -\varepsilon)y}.$ It is also clear that
\begin{multline*}
	{\sf P}\left(|B_1|>t_1\right)\leqslant\binom{n}{t_1}q_1^{t_1}\leqslant (nq_1)^{t_1}\leqslant e^{-\frac{\varepsilon^2}{25}y(1-o(1))\ln^2\ln n},\\
	{\sf P}\left(|B_2|>t_2\right)\leqslant\binom{n}{t_2}q_2^{t_2}\leqslant (nq_2)^{t_2}\leqslant e^{-\left(\frac{1}{\ln 2}+3\varepsilon\right)\left(\ln 2 - \varepsilon\right)y(1-o(1))\ln\ln n}<e^{-\left(1+\frac{\varepsilon}{2}\right)y\ln\ln n}.
\end{multline*}
Expose all edges between $Y$ and $[n]\setminus Y.$ Assume that, for every $v\in [n]\setminus Y,$ we have 
$$
N_0(v,Y)\geqslant (\ln n)^{2+2\varepsilon}\sqrt{n},\qquad
b_1:=|B_1|\leqslant t_1,\qquad b_2:=|B_2|\leqslant t_2.
$$ Put $$
m := \ceil[\big]{(\ln n)^{1+\varepsilon}\sqrt{n}}.
$$
For convenience, set $B_2=\{1,\ldots,b_2\}.$ For every $v\in B_2,$ choose $Z_v\subset Y$ such that $|Z_v|=m,$ $Z_v\subset N_0(v,Y)$ and $Z_1,\ldots,Z_{b_2}$ are disjoint.

Set $Y' = Y\setminus \left(Z_1\sqcup\ldots\sqcup Z_{b_2}\right),$ $y':=|Y'|=y-b_2m.$ For every $v\in [n]\setminus\left(Y\sqcup B_1\sqcup B_2\right),$
\begin{equation}\label{eq1}
	\left|N_0\left(v,Y'\right)\right|\geqslant\left(\frac{1}{2}-\frac{\varepsilon}{5}\frac{\ln\ln n}{\ln n}\right)y-b_2m=\left(\frac{1}{2}-\frac{\varepsilon}{5}\frac{\ln\ln n}{\ln n}\right)y'-O\left(t_2m\right).
\end{equation}
For every $v\in B_1,$
\begin{equation}\label{eq2}
	\left|N_0(v,Y')\right|>\alpha y - b_2m = \alpha y' - O\left(t_2m\right).
\end{equation}

Choose $k-b_2$ disjoint sets $Z_{b_2+1},\ldots,Z_k$ in $Y'$ of size $m$ uniformly at random. Let us show that whp every vertex of $[n]\setminus\left(Y\sqcup B_1\sqcup B_2\right)$ has at least roughly $\frac{1}{2}|Z_i|$ non-neighbors in $Z_i$ for every $i\in[k]\setminus[b_2].$ Let $v\in[n]\setminus\left(Y\sqcup B_1\sqcup B_2\right).$ If we choose sets $Z_{b_2+1},\ldots,Z_k$ sequentially then, for each $i\in[k]\setminus[b_2],$ the conditional distribution of the number of non-neighbors of $v$ in $Z_i$ given $Z_{b_2+1},\ldots,Z_{i-1}$ such that $\left||N_0(v,Z_j)|-|N_0(v,Y')|\cdot\frac{m}{y'}\right|\leqslant m^{3/4},$ is hypergeometric with expectation 
\begin{equation}\label{eq3}
	\lambda_i\geqslant \frac{m}{y'}|N_0(v,Y')| -\frac{(i-b_2-1)m^{7/4}}{y'-(i-b_2-1)m}.
\end{equation} By Hoeffding's bound \cite[Theorem 2.10]{JLR}, the probability that this number is outside $\left[\lambda_i - m^{3/4}(1-\varepsilon),\lambda_i+m^{3/4}(1-\varepsilon)\right]$ is at most $e^{-\Omega\left(\sqrt{m}\right)}.$ By the union bound, the probability of existence of a vertex $v\in[n]\setminus\left(Y\sqcup B_1\sqcup B_2\right)$ and $i\in[k]\setminus[b_2]$ such that $v$ has less than $\lambda_i - m^{3/4}$ non-neighbors in $Z_i$ is at most $nke^{-\Omega\left(\sqrt{m}\right)}\to 0$ as $n\to\infty.$

Similarly, by the union bound and by \cite[Theorem 2.10]{JLR}, the probability that there exists a vertex $v\in B_1$ and $i\in[k]\setminus[b_2]$ such that $v$ has less than $\frac{m}{y'}|N_0(v,Y')|-m^{3/4}$ non-neighbors in $Z_i$ is at most $b_1ke^{-\Omega(\sqrt{m})}\to 0$ as $n\to\infty.$

So, by $\eqref{eq1}-\eqref{eq3},$ for $n$ large enough, we may choose disjoint subsets $Z_{b_2+1},\ldots,Z_k\subset Y'$ of size $m$ s.t. each $v\in[n]\setminus \left(Y\sqcup B_1\sqcup B_2\right)$ has at least $\left(\frac{1}{2}-\frac{\varepsilon\ln\ln n}{5\ln n}\right)m - O\left(\frac{t_2m^2}{y}\right)$ non-neighbors in each $Z_i$ and each $v\in B_1$ has at least $\alpha m - O\left(\frac{t_2m^2}{y}\right)$ non-neighbors in each $Z_i,$ $i\in[k]\setminus[b_2].$ Let $Z_{b_2+1},\ldots,Z_{k}$ be such sets.

Denote by $\mathcal{F}$ the family of all subsets of size $k$ of $Y$ that contain one element in each set $Z_i,$ $i\in[k],$ and contain at least one non-neighbor of each vertex $v\in [n]\setminus Y.$ By the definition of $Z_i,$ $i\in[b_2],$ any element of $\mathcal{F}$ contains at least one non-neighbor of every $v\in B_2.$ Consider the uniform distribution on the set of all $k$-sets in $Y$ that contain one element in each $Z_i,$ $i\in[k].$ The probability of $\overline{\mathcal{F}}$ equals the probability that there exists a vertex in $[n]\setminus \left(Y\sqcup B_2\right),$ connected to all vertices of the random $k$-set. Then, by the union bound,
\begin{multline*}
	1-\frac{|\mathcal{F}|}{m^k}\leqslant n\left( \frac{1}{2}+\frac{\varepsilon}{5}\frac{\ln\ln n}{\ln n} + O\left(\frac{t_2m}{y}\right)\right)^{k-b_2}+b_1\left(1-\alpha+O\left(\frac{t_2m}{y}\right)\right)^{k-b_2}\leqslant\\ ne^{\left(\log n +\left(\frac{1}{\ln 2} + 4\varepsilon\right)\ln\ln n - \left(\frac{1}{\ln 2} + 3\varepsilon\right)\ln\ln n\right)\left(-\ln 2 +\frac{2\varepsilon}{5}\frac{\ln\ln n}{\ln n}(1+o(1))\right)}+\\
	\ln^2n\cdot e^{\left(\log n +O\left(\ln\ln n\right)\right)\left(-\ln\frac{1}{1-\alpha}+o\left(1\right)\right)}=\\
	e^{-\varepsilon\left(\ln 2 - \frac{2}{5\ln 2}-o(1)\right)\ln\ln n}+o(1)\to 0
\end{multline*}
as $n\to\infty.$

Therefore 
$|\mathcal{F}|\geqslant m^k(1-o(1)).$

It remains to prove that, if we expose the edges inside $Y,$ then with probability $1-O\left(e^{-\left(1+\frac{\varepsilon}{2}\right)y\ln\ln n}\right)$ at least one set from $\mathcal{F}$ induces a clique.

For each member $K$ of $\mathcal{F}$ denote by $\mathbb{I}_K$ the indicator random variable equal to $1$ iff $K$ induces a clique in $G(n,\frac{1}{2}).$ Define $\eta:=\sum\limits_{K\in\mathcal{F}}\mathbb{I}_K,$ $\mu:={\sf E}\eta = |\mathcal{F}|2^{-\binom{k}{2}}  =(1-o(1))m^k2^{-\binom{k}{2}}.$ Define $\Delta := \sum\limits_{K,\;K'}{\sf P}\left(\mathbb{I}_K=\mathbb{I}_{K'}=1 \right),$ where the summation is over all ordered pairs $K,$ $K'$ of members of $\mathcal{F}$ that satisfy $2\leqslant|K\cap K'|\leqslant k-1.$ Denote by $\Delta_i$ the contribution of pairs $K,K'\in \mathcal{F}$ for which $|K\cap K'|=i.$ Then
\begin{multline*}
	\frac{\Delta}{\mu^2}=\sum_{i=2}^{k-1}\frac{\Delta_i}{\mu^2}\leqslant \sum_{i=2}^{k-1}\frac{|\mathcal{F}|}{\mu^2}2^{-\binom{k}{2}}\binom{k}{i}(m-1)^{k-i}2^{-\binom{k}{2}+\binom{i}{2}}\leqslant\\ \leqslant \sum_{i=2}^{k-1}\frac{m^{k-i}2^{\binom{i}{2}}}{|\mathcal{F}|}\binom{k}{i}\leqslant\left( \frac{k^2}{m^2}+\sum_{i=3}^{k-1}\left(\frac{2^{\frac{i}{2}}ke}{im} \right)^i\right) (1+o(1)).
\end{multline*}
Observe that the sign of
\begin{equation*}
	\frac{\partial\ln\left[ \left(\frac{2^{\frac{i}{2}}ke}{im} \right)^i\right] }{\partial i} = i\ln 2 - \ln \frac{m}{k}-\ln i
\end{equation*}
changes on $[3,k-1]$ only once from negative to positive, which implies that  $\left(2^{i/2}ke/im\right)^i$ as a function of $i$ attains its maximum at one of the endpoints. If $i = 3,$ then 
\begin{equation*}
	\left(\frac{2^{\frac{i}{2}}ke}{im}\right)^i = O\left(\frac{k^3}{m^3}\right).
\end{equation*}
If $i = k-1,$ then 
\begin{equation}\label{janson}
	\left(\frac{2^{\frac{i}{2}}ke}{im}\right)^i\leqslant e\left(\frac{2^{\frac{k-1}{2}}e}{m}\right)^{k-1} \leqslant e\left(\frac{e\sqrt{n}\left(\ln n\right)^{\frac{1}{2\ln 2}+2\varepsilon}}{m}\right)^{k-1} \leqslant e\left(\frac{e}{(\ln n)^{\varepsilon}}\right)^{k-1}(1+o(1)).
\end{equation}
Hence, $\frac{\Delta}{\mu^2}\leqslant \frac{k^2}{m^2}(1+o(1)).$ Also, $\frac{1}{\mu} = \frac{2^{\binom{k}{2}}}{m^k}(1+o(1)) = o\left(\frac{k^2}{m^2}\right).$
By Janson's inequality \cite[Theorem 2.18]{JLR},
\begin{equation*}
	{\sf P}\left(\eta=0 \right) \leqslant e^{-\frac{\mu^2}{\mu + \Delta}}\leqslant e^{-\frac{(1-o(1))m^2}{2k^2}}=e^{-\Omega\left(n(\ln n)^{4\varepsilon}\right)}.
\end{equation*}
This completes the proof of Lemma \ref{lemma}.

\section{Conclusion}\label{conclusion}
We have shown that the clique chromatic number $\chi_c\left(G\left(n,\frac{1}{2}\right) \right)$ is, whp, $\frac{\log n}{2} - \Theta\left(\ln\ln n\right).$ The same proof can be applied to the random graph $G(n,p)$ with any constant probability $\frac{1}{2}\leqslant p<1.$ So, whp $\chi_c\left(G\left(n,p\right) \right)=\frac{\log_{1/(1-p)}n}{2} - \Theta\left(\ln\ln n\right),$ $\frac{1}{2}\leqslant p<1.$ Observe that the bounds in \eqref{janson} can not be directly transferred to the values of $p$ less than $\frac{1}{2}.$ So the problem of finding the asymptotics of $\chi_c\left(G\left(n,p\right) \right)$ for constant $p\in\left(0, \frac{1}{2}\right)$ remains open.

Notice that our arguments also imply that all clique colorings in not too many colors are close to the ``greedy coloring''. Formally, for every $\varepsilon>0,$ whp for every coloring of the vertex set of $G\left(n,p \right)$ in $s\leqslant\log_{\frac{1}{1-p}} n$ colors, there exists an ordering $Y_1,\ldots,Y_s$ of the color classes such that, for every $j<\frac{\log_{1/(1-p)} n}{2} - O\left(\ln\ln n\right),$ there exist vertices $v_1,\ldots,v_j$ such that $|\left( Y_1\cup\ldots\cup Y_j\right) \setminus \left( N(v_1)\cup\ldots\cup N(v_j)\right)|\leqslant \sqrt{n}\left(\ln n\right)^{\varepsilon},$ where $N(v)$ is the neighborhood of $v$ in $G(n,p).$ Also note that the upper bound that we prove in Section \ref{section_upper_bound} is obtained by the following greedy procedure: for every $i=1,\ldots,s,$ assign the color $i$ to all neighbors of vertex $i$ that did not receive colors $1,\ldots,i-1.$ This algorithm yields an asymptotically tight (up to an additive term $O(\ln\ln n)$) bound. This differs from the outcomes of greedy algorithms for finding the ordinary chromatic number: the greedy coloring requires roughly twice the number of colors in a proper coloring of $G(n,p),$ see \cite{GM}.

\section*{Appendix: proof of Claim \ref{lemma_upper}}
Let $X_k$ be the number of cliques of size $k=k(n),$ $k\in\mathbb{N},$ in $G(n,1/2)\vert_{\left[m\right]}$ such that any vertex outside $[m]$ has a non-neighbor inside the clique. We have
\begin{equation*}
	{\sf E}X_k = \binom{m}{k}\left(\frac{1}{2}\right)^{\binom{k}{2}}\left(1 - \frac{1}{2^k}\right)^{n-m}\leqslant e^{k\ln m - k\ln k + k - \binom{k}{2}\ln 2 - \frac{n-m}{2^k}}.
\end{equation*}
Let $f(k) = k\ln m - k\ln k + k - \binom{k}{2}\ln 2 - \frac{n-m}{2^k}.$
Clearly,
\begin{equation*}
	\frac{\partial f}{\partial k} = \ln m - \ln k - k\ln 2 + \frac{\ln 2}{2} + \frac{n-m}{2^k}\ln 2
\end{equation*}
decreases and $k_0 = \log n - \log\ln \frac{n}{m} + O(1)$ is the only solution of $\partial f/\partial k = 0.$ Therefore, $f(k)\leqslant f(k_0).$ We obtain
\begin{multline*}
	f(k_0) = \log n\cdot\ln m - \log\ln\frac{n}{m}\cdot\ln m - \log n\cdot\ln\ln n -\\ \frac{1}{2}\log^2 n\cdot\ln 2 + \log n\cdot\log\ln\frac{n}{m}\cdot\ln 2 + O\left(\ln n\right).
\end{multline*}
Observe that $f(k_0)$ is an increasing function of $m.$ Hence, we have
\begin{multline*}
	f(k_0)\leqslant \frac{1}{2}\log n\cdot\ln n + \left(\frac{1}{2}-\varepsilon\right)\log n\cdot \ln\ln n - \frac{1}{2}\log \ln n\cdot\ln n - \\ \log n\cdot\ln\ln n - \frac{1}{2}\log n \cdot\ln n + \log n\cdot\log\ln n\cdot\ln 2 + O(\ln n) =\\
	-\varepsilon(1+o(1))\log n\cdot\ln\ln n.
\end{multline*}

Finally, we obtain that the probability of the existence of two positive integers $m\leqslant \sqrt{n}(\ln n)^{1/2 - \varepsilon}+\sqrt[4]{n}\left(\ln n\right)^{5/4 - \varepsilon / 2}$ and $k\leqslant m$ such that $X_k>0$ is at most
\begin{equation*}
	\left(\sqrt{n\ln n}\right)^2e^{-\varepsilon(1+o(1))\log n\cdot\ln\ln n}\to 0
\end{equation*}
as $n\to\infty.$
\end{document}